\documentclass[12pt]{amsart}
\usepackage{url}
\usepackage{tikz}

\usetikzlibrary{arrows,shapes,automata,backgrounds,decorations,petri,positioning,patterns}

\numberwithin{equation}{section}

\newtheorem*{unnumthm}{Theorem}

\theoremstyle{definition}

\begin{document}

\bibliographystyle{siam}

\title{Fingerprint Databases for Theorems}
\author[Sara C. Billey]{Sara C. Billey$^1$}
\thanks{$^1$Research partially supported by a grant from the National Science Foundation DMS-1101017.}
\address{University of Washington, Mathematics Department, Box 354350, Seattle, WA 98195 }
\email{billey@uw.edu} 
\author[Bridget E. Tenner]{Bridget E. Tenner$^2$}
\thanks{$^2$Research partially supported by a DePaul University Competitive Research Leave grant.} 
\address{Department of Mathematical Sciences, DePaul University, Chicago, IL 60614}
\email{bridget@math.depaul.edu}

\maketitle

\begin{quote}
``\emph{Fingerprint}, in the anatomical sense, is a mark made by the
pattern of ridges on the pad of a human finger.  The term has been
extended by metaphor to anything that can uniquely distinguish a
person or object from another.''  \cite{wiki.fingerprint}.
\end{quote}

Suppose that $M$ is a mathematician, and that $M$ has just proved
theorem $T$.  How is $M$ to know if her result is truly new, or if $T$
(or perhaps some equivalent reformulation of $T$) already exists in
the literature?  In general, answering this question is a nontrivial
feat, and mistakes sometimes occur.

Certain mathematical results
have canonical representations, or \emph{fingerprints}, and some
families of fingerprints have been collected into searchable
\emph{databases}.  If $T$ is such a theorem, then $M$'s search will be
greatly simplified.  Note that the searchable nature of a database is important here.  An analogue of
``alphabetical order'' does not exist for all structures, and so it is
important that $M$ is able to query the fingerprint of $T$, instead of
needing to browse through all existing catalogued results.

A revolutionary mathematical tool appeared online in 1996 --- Neil
Sloane's collection of integer sequences, along with mathematical
interpretations of the numbers, formulas for generating them, computer
code, references, and relevant links.  This was the On-Line
Encyclopedia of Integer Sequences (OEIS) \cite{oeis}, originally
hosted on Sloane's website at AT\&T Labs.  Anyone with access to the
internet could peruse the database, and anyone could submit a sequence
or supplemental data to the database.  All for free.  Thanks to
Sloane's tireless efforts and a worldwide community of contributors,
the collection has grown to well over 200,000 sequences to date,
drawing results from all areas of mathematics.  Each sequence in the
OEIS acts as a fingerprint for an associated theorem.  While the
fingerprints in the OEIS have a specific input structure, the
sequences can arise in many contexts, including arrays of data,
coefficients of polynomials, enumeration problems, subway stops, and
so on.  The OEIS itself is the database for these fingerprints.  The
impact on research is clearly established by over 3000 articles to
date citing the OEIS \cite{oeis.citations}.

Fingerprinting has made an impact in many scientific fields.  For
example, fingerprinting documents is crucial in computer science for reducing duplication
in web search results, isotopic fingerprints are used in fields
ranging from chemistry to archaeology, and there is of course
extensive use of fingerprinting in forensic science.

There are other families of mathematical results that have their own
identifying fingerprints, not in the form of integer sequences.
Searchable catalogues are already in use for some of these families,
while no such directories yet exist for others.  The aim of this
article is to give these resources greater exposure, and also to
encourage the community to create and support new fingerprint
databases for other mathematical structures.  Note that what we
propose are not simply enhanced digital mathematical resources.
Rather, a fingerprint database of theorems should be a searchable,
collaborative database of citable mathematical results indexed by small,
language-independent and canonical data.

Every day new tools for searching the scientific literature are
established.  To be clear, this article will be out of date the moment
it is published.  In fact, active research at the intersection of
mathematics, computer science, and linguistics is devoted to
organizing mathematics into more searchable formats, including the
Mathematical Knowledge Management and Intelligent Computer Mathematics
conferences.  An example outside of mathematics is biomedical natural language
processing, known as BioNLP \cite{bionlp}.  We expect that one
day, natural language processing will be applicable to theorems and
will significantly facilitate $M$'s search through the literature.
The question is, what can we do until then?

Rod Brooks and his group at MIT used the phrase
\begin{center}
``fast, cheap, and out of control''
\end{center}
to describe an emphasis on building small, cheap, and redundant robots instead of overly
complex single machines \cite{brooks}.  We suggest that a similar
approach to fingerprinting theorems can make a big impact in the near
future, while more finessed tools are being developed in the background.  It is better to start
a theorem collection now --- with an imperfect but efficient fingerprint --- than to waste time awaiting an epiphany about the perfect mechanism for encoding this data.

In a sense, we are proposing a new line of research for mathematicians to address: what are the fingerprintable theorems within each discipline of mathematics, and what might those fingerprints look like?

\section*{Known results can be hard to find}

Theorems are usually written in human-readable language.  They employ
specialized vocabulary, functions, and layers of hypotheses and
implications.  A theorem in one branch of mathematics can resurface in
another context, and the two statements may bear little superficial
resemblance to each other.  Search engines can help uncover a result
if it is accessible online and there is a name associated to it, such
as for a solution to a famous conjecture in which case the name, or
names, would be the fingerprint.  For example, one can easily ask a
search engine for information about Fermat's Last Theorem, which would
lead $M$ to discover that her result $T$ was already proved by Wiles
\cite{wiles}.  

Formulas are prevalent in mathematical research, but are inherently difficult to query.  For example, $M$ would have to make decisions about notation, variable names, and formatting.  Moreover, even if search engines did have a good mechanism for querying formulas, it might not be especially useful --- a given formula can often be stated in a variety of ways.  For example, the following basic trigonometric identities are equivalent:
$$\begin{array}{c@{\ +\ }c@{\ =\ }c@{}c}
\sin^2\theta & \cos^2\theta &1,&\\
2\tan^2\theta & 2 & 2\sec^2\theta, &\text{ and}\\
3 & 3\cot^2\theta & 3\csc^2\theta.&
\end{array}$$
If our mathematician $M$ has discovered a new statement of an existing formula, a search engine might have difficulty detecting that her result is equivalent to the known one.

There have been many ideas put forth for improving the search tools
for formulas in the literature.  In fact, search tools themselves can
contribute to mathematical results.  Notably, G\"odel invented a
numerical encoding of formulas as a step toward proving his famous
Incompleteness Theorem \cite{godel}.  However, the procedure is not
unique and it is certainly not efficient.  For example, the G\"odel
number of the formula ``0 = 0'' is $2^6\times 3^5 \times 5^6 = 243,000,000$.  More
recently, Borwein and Macklem address the question of how best to add
hyperlinks to electronically available textbooks \cite{Borwein-Macklem}.

Of course, $M$'s search through existing literature for any hint of
her theorem $T$ would have been much harder prior to the internet.
There are examples throughout mathematical history of theorems having
been discovered, and subsequently rediscovered independently ---
sometimes over and over again.  For example, the characterization of
higher-dimensional regular polytopes, attributed to Schl\"afli, had been 
recovered at least nine other times by the end of the nineteenth
century \cite{coxeter}.

\section*{Benefits of a good fingerprint database}

We wish to proselytize the accumulation of theorem fingerprints into databases.  We urge the reader to become a collector.  A connoisseur, even! 

First, though, we must explain how the OEIS encodes \emph{theorems} --- after all, its primary purpose is to collect and catalogue integer sequences.  In fact, the theorems can be found within the architecture of this database --- namely by the inclusion of other fields associated to each sequence such as ``name,'' ``comments,'' ``formula,'' and so on.

If our mathematician $M$ is going to make use of the OEIS, it is because she has encountered a sequence of integer data within her work.  Then $M$ runs a query against the OEIS, using her data.  Even a relatively small subsequence --- perhaps just two numbers --- can sometimes determine a unique entry in the OEIS.  The responses from $M$'s search enable her to connect her data to known literature, to find formulas, to make conjectures, and so on.

For example, if $M$ enters $0,1,1,2,3,5,8,13,21$ into the OEIS, the first option it returns is sequence A000045, the Fibonacci numbers.  Two of the comments for this entry are
\begin{itemize}
\item $F(n+2) = $ number of subsets of $\{1,2,\ldots,n\}$ that contain no consecutive integers, and
\item $F(n+1) = $ number of tilings of a $2\times n$ rectangle by $2\times 1$ dominoes.
\end{itemize}
Thus this entry encodes a variety of results, including the following.
\begin{quote}\vspace{-.25in}
\begin{unnumthm}[{\cite[A000045]{oeis}}]
The subsets of $\{1,2,\ldots,n\}$ containing no consecutive integers are in bijection with the tilings of a $2\times (n+1)$ rectangle by $2\times 1$ dominoes, and these are each enumerated by the $(n+2)$nd Fibonacci number.
\end{unnumthm}
\end{quote}
In this way, each entry in the OEIS chronicles a mathematical theorem, and the integer sequence associated to the entry is that theorem's fingerprint.  The OEIS is arguably the most established fingerprint database for theorems to date.

\section*{Other fingerprint databases}

Depending on the structure of theorem $T$, the OEIS is not the only
tool of its kind available to the curious $M$.  We will describe some
of the fingerprint databases for theorems that already exist in this
section.  These databases augment the classical approach to finding
theorems in the literature, including books, journals, MathSciNet, the arXiv,
and the World Digital Mathematics Library.

\subsection*{Permutation patterns}
The Database of Permutation Pattern Avoidance (DPPA) \cite{dppa} contains collections of permutations --- thought of as patterns --- whose avoidance exactly characterize particular phenomena.  The second author started this database in 2005, and it has grown to more than 40 sets of patterns so far.  In addition to the patterns themselves, each entry in the DPPA includes the phenomenon (or phenomena) being characterized, references to existing literature, and a link to the OEIS whenever possible.  The DPPA is searchable both by permutation (pattern) and by keyword.

For example, if the theorem $T$ involves permutations avoiding the two patterns $3412$ and $4231$, then the DPPA would have directed $M$ to entry P0005, for the set $\{3412,4231\}$.  The two descriptions for this entry are
\begin{itemize}
\item permutations with rank symmetric order ideals in the Bruhat order, and
\item permutations indexing smooth Schubert varieties,
\end{itemize}
as described in \cite{carrell, lakshmibai sandhya}.

Each entry of the DPPA represents a characterization theorem.  The theorem for the entry just described would be as follows.
\begin{quote}\vspace{-.25in}
\begin{unnumthm}[{\cite[P0005]{dppa}}]
The permutations with rank symmetric order ideals in the Bruhat order are exactly those that index smooth Schubert varieties, and they are precisely the permutations that avoid the patterns $3412$ and $4231$.
\end{unnumthm}
\end{quote}
The fingerprint for each DPPA theorem is its associated set of patterns, and the DPPA itself is the database for these fingerprints.

\subsection*{FindStat} FindStat \cite{findstat} is a database of
statistics on combinatorial objects.  It was created in 2011 by Berg
and Stump, and currently catalogues over $50$ statistics.  If $M$ has
obtained some data about one of these objects, then she could enter
her data into FindStat, and it would tell her if this particular
statistic is included in the database.  If so, FindStat would identify
the standard vocabulary used for that statistic.  This would equip $M$
with searchable terminology, allowing her to discover any relevant
existing literature.

\subsection*{Hypergeometric series} Every hypergeometric series can be
written in a canonical form, and this form serves as the fingerprint
for these objects.  It has long been common to store identities for
these series in tables, listed in a given order by these canonical
forms.  For example, Bailey published such a collection in 1935
\cite{bailey}.  Perhaps this book is the original fingerprint database
for theorems.

The modern approach has taken research in hypergeometric identities
one step further.  The WZ method for finding identities involving
hypergeometric series has been described in the book \textit{A=B} by
Petkov\v{s}ek, Wilf, and Zeilberger \cite{petkovsek wilf zeilberger}.
Using these algorithms, one can determine definitively if a
hypergeometric series has a closed form or not.  If there is a closed
form, then the WZ method will produce it, given enough computational
time and memory.  Furthermore, this procedure will give a proof
certificate that can be used to check the identity.  Many new
identities and new proofs of known identities have been found using
the WZ method, for example \cite{ekhad zeilberger}.  What this
resource currently lacks is a way to connect results to existing
literature, pointing our mathematician $M$ to what is already known
about each identity.  

The NIST Digital Library of Mathematical Functions (DLMF) also
includes many hypergeometric identities indexed by canonical form and
some references.  We should point out, however, that neither the WZ
method nor the DLMF form a fingerprint database for theorems
themselves in their current form.  Perhaps there could be a
collaborative effort to catalogue all known hypergeometric identities
with extensive references, and entries searchable by their canonical
forms.  If so, all new identities found by the WZ method could include
their proof certificate as a comment.  This could provide a useful
place to ``publish'' proof certificates.

\section*{Constructing a fingerprint database is not always easy}

An important asset of the OEIS, the DPPA, FindStat, and the WZ method
is that the fingerprints they use are \emph{language independent}.
More precisely, their input is entirely numerical and canonical ---
free from specialized vocabulary.  This seems to be a necessary
feature of a good fingerprint database for theorems.

Another desirable feature of a productive fingerprint database is that
it should \emph{reference} existing literature whenever possible.
Cross-references within a single database and between different
databases can only enhance the state of knowledge.  Features like
computer code and external links can be highly beneficial when
relevant.  For example, any integer sequence associated with a theorem
in a new fingerprint database should reference the relevant OEIS
entry.

Because mathematics is so broad and develops so quickly, a fingerprint
database for theorems should be \emph{collaborative} --- publicly
available and  welcoming additions from anyone subject to editorial
standards.  The Wikipedia model for an open database is a highly
successful model of this idea.  However, one does not need to learn
MediaWiki before starting a collection of theorem fingerprints;
rather, one could simply ask for new database entries to be submitted
in some kind of standard format which can easily be added to the
database.

Finally, it is most convenient for the fingerprint to be encoded in a
\emph{small} amount of data.  There is a natural conflict between
keeping fingerprints small and uniquely identifying each
object in the database.  Certainly some compromises to one or both of
these might be necessary.  An efficient fingerprint encryption may be permitted to
return some false positives, but it should never return a false negative.
The possibility of false positives is all the more reason for additional fields
within the database entries, to distinguish the true from the false
positives.  For example, querying the first nine Fibonacci numbers will return many false positives in the OEIS, but $M$ can weed through them by reading through their full records.

There are certainly some challenges to creating a fingerprint
database for theorems.  These include identifying the right data structure as the
fingerprint, determining a canonical format, addressing structures that
have no obvious numerical data, and compactly encoding
a given fingerprint.  We hope these obstacles will not be too daunting, though, because an imperfect resource is still better than no resource at all.  Two examples are given below.

\subsection*{Example: fingerprinting graphs}

Theorems about finite graphs deserve a fingerprint database.  There exist
numerous classification theorems in graph theory that equate graph
containment with important properties.  One of the monumental results of the twentieth
century is the Graph Minor Theorem by Robertson and Seymour
\cite{robertson seymour}:
\begin{quote}
Any family $\mathcal{F}$ of graphs that is closed under taking minors
can be characterized as the set of all graphs whose minors avoid a
finite list $L(\mathcal{F})$.
\end{quote}
This result certainly suggests that graphs can fingerprint theorems.  The Wagner formulation of Kuratowski's Theorem is an example of this situation \cite{kuratowski, wagner}:  
\begin{quote}
A simple graph $G$ is planar if and only if $G$ has no minor
isomorphic to the graphs known as $K_{3,3}$ and $K_5$.
\end{quote}
Graphs arise as classification tools in many fields of
mathematics, including Hales's proof of Kepler's Conjecture
\cite{hales} and the classification of finite Coxeter groups
\cite[Chapter 2 and Section 6.4]{hcox}.  

One could enumerate the results of a graph theorem, say by counting the graphs of each size possessing a certain property.  The resulting sequence could be an entry in the OEIS.  However, a graph theorem database would still be relevant because it could track more specific graph properties through further refinement and cataloguing.  Moreover, and perhaps more persuasively, counting graphs is not an easy computational problem, so this partial enumerative fingerprint would not uniquely identify the appropriate entry in the OEIS.  For example, the linklessly embeddable graphs in Euclidean space are characterized by avoiding the Petersen family of graphs, which include seven graphs having between six and ten vertices each.  It is computationally infeasible to compute the number of linklessly embeddable graphs on six, seven, eight, nine, and ten vertices, which would be the first few times at which this sequence differs from the sequence enumerating all graphs.

There currently exist many online resources for graph data, such as House of Graphs \cite{hog} and the tools listed at \cite{wiki.graph.database}.  However, none of these resources are databases of theorems (at present).  It is inherently difficult to fingerprint graph theorems using
searchable, canonical, and concise numerical data.  In particular,
there is not an obvious choice for the best way to fingerprint a
graph.

The adjacency matrix of a graph describes the graph uniquely in numerical data.  Often in graph theory, a classification theorem depends only on
isomorphism classes.  This could pose a problem if the fingerprint of
a graph is its adjacency matrix because isomorphic graphs can have
different adjacency matrices.  For example, the graph with two
adjacent vertices and one isolated vertex could be represented by
any of
$$\begin{bmatrix}
0 & 1 & 0\\
1 & 0 & 0\\
0 & 0 & 0
\end{bmatrix},
\hspace{.25in}
\begin{bmatrix}
0 & 0 & 1\\
0 & 0 & 0\\
1 & 0 & 0
\end{bmatrix},
\hspace{.25in}
\text{and}
\hspace{.25in}
\begin{bmatrix}
0 & 0 & 0\\
0 & 0 & 1\\
0 & 1 & 0
\end{bmatrix}.$$
We can, of course, handle this difficulty by choosing a canonical
representative in each isomorphism class, such as the adjacency matrix
whose row reading word is smallest in lexicographic order.  However,
finding such a canonical adjacency matrix is
no easy task: there is no known polynomial time algorithm for testing graph
isomorphism.  In fact, it is an open question whether the graph isomorphism
problem is NP-complete.

Degree sequences are an attractive choice for fingerprints because they are much easier to encode than adjacency matrices.  If one were to fingerprint graph families by lists
of degree sequences written in lexicographic order, then $K_{3,3}$ and
$K_{5}$ would be encoded as the list $[[3,3,3,3,3,3], [4,4,4,4,4]]$.  Querying this list in a database of graph theorems, the mathematician $M$ would learn that these two graphs are
related to planar graphs via Kuratowski's Theorem.

On the other hand, a degree sequence does not determine a unique graph.  For example, both
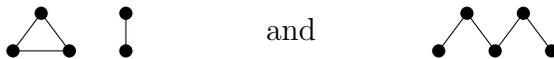
\begin{center}
\begin{minipage}{.75in}
\begin{tikzpicture}[scale=.5]
\draw (0,0) -- (1.5,0) -- (.75,1) -- (0,0);
\draw (3,0) -- (3,1);
\foreach \x in {(0,0),(1.5,0),(.75,1),(3,0),(3,1)} {\fill[black] \x circle (5pt);}
\end{tikzpicture}
\end{minipage}
\hspace{.5in}
\text{and}
\hspace{.5in}
\begin{minipage}{.75in}
\begin{tikzpicture}[scale=.5]
\draw (0,0) -- (.75,1) -- (1.5,0) -- (2.25,1) -- (3,0);
\foreach \x in {(0,0),(1.5,0),(.75,1),(3,0),(2.25,1)} {\fill[black] \x circle (5pt);}
\end{tikzpicture}
\end{minipage}
\end{center}
have degree sequence $[2,2,2,1,1]$.  

Theorems about specific graphs or families of graphs may be rare enough that the compromises one makes when fingerprinting by degree sequences might not result in too many false positives.  Indeed, as we have said before, it is better to have a collection
of theorems with an imperfect fingerprint, than to have no collection
at all!

\subsection*{Example: finite groups}
The finite simple groups have been completely classified
\cite{wilson}.  These groups fall into six families, and the title for
each group is given by a combination of letters and numbers.  For
example, one group is denoted $^3D_4(q^3)$.  These groups, and various
details about them, are collected in the ATLAS of Finite Group
Representations \cite{atlas}.  To date, this resource includes more
than $5000$ representations of more than $700$ groups.

The current implementation of the ATLAS does not allow users to search
the database by numerical invariants of the groups, thus it is not a
fingerprint database as we have defined it.  To find the details of a
group, one must know its title or something about where it fits into the
classification.

To make the ATLAS into a fingerprint database, one would have to add a
feature where groups could be detected by some numerical
invariant(s).  For example, an additional search box could be added to
the main webpage to access the database by entering the order of a
group.  Then the order would act as the fingerprint.  There are groups
of the same order already in the database, but perhaps the number of
coincidences is small enough that a user could prune the results via
the many other entries available.  Additional invariants
might also be used to refine the search.

\section*{What should happen next}

We believe that many families of theorems can be fingerprinted ---
some identified by obvious data structures, others perhaps by less
obvious ones.  We encourage everyone in the mathematical community to
look in their own work for results that can be identified by some
form of compact data.  In fact, any structure that has a canonical
parameterization merits this attention.  Additionally, a long term
benefit of having these databases is that structures amenable to
fingerprinting may also be amenable to computer proof verification systems
and learning algorithms, as with the Four Color Theorem
\cite{gonthier,robertson sanders seymour thomas} and permutation
patterns \cite{MU}.

Clever insight, beyond what is currently common practice, might be necessary to find an appropriate fingerprint.  In fact, the need to find theorem fingerprints can drive future research.

Many disciplines of mathematics would benefit from the greater context of a theorem
database.  The accessibility of mathematical research in the last few decades has
flourished.  In the past few years alone, we have seen substantial
growth as measured in mathematics articles posted on the arXiv, increasing from
4654 articles in 2002 to 24176 articles in 2012 \cite{arxiv}.  With
this level of productivity, fingerprint databases are even more
valuable.  These resources --- both the ones that currently exist and those
that we hope the readers will create --- enhance experimental
mathematics, help researchers make unexpected connections between
areas, and even improve
the refereeing process.  We encourage everyone to follow Neil Sloane's lead and to take up such a collection.\\

Hats off to Neil!

\section*{Acknowledgments} First and foremost, we want to thank the
OEIS and all of its contributors, with special thanks to Neil Sloane.
We also thank all the contributors to the other resources we have
referenced and used in our own work.  We would like to thank Chris
Berg, Jon Borwein, Neil Calkin, Chris Godsil, Ron Graham, Ursula
Martin, Brendan Pawlowski, Christian Stump, Lucy Vanderwende, Paul Viola, and Doron
Zeilberger for helpful discussions while preparing this article.  We
thank the organizers of the ICERM workshop on Reproducibility in
Computational and Experimental Mathematics for presenting a chance to
discuss these ideas with a broad community.  Finally, the first author
thanks Rod Brooks for the opportunity to work in his lab as an
undergraduate at the height of the ``fast, cheap, and out of control''
revolution in robotics.

\end{document}